\documentclass[11pt]{article}

\usepackage{amsfonts}
\usepackage{amsmath}
\usepackage{amssymb}
\usepackage{hyperref}

\newtheorem{theorem}{Theorem}
\newtheorem{corollary}{Corollary}

\newtheorem{lemma}{Lemma}
\newtheorem{proposition}{Proposition}
\newtheorem{remark}{Remark}
\newenvironment{proof}[1][Proof]{\noindent\textbf{#1.} }{\ \rule{0.5em}{0.5em}}

\def\v{\nu}
\def\I{\infty}

\begin{document}

\title{\textbf{On the $q$-Bessel Fourier transform}}
\author{ Lazhar Dhaouadi \thanks{
IPEIB, 7021 Zarzouna, Bizerte,Tunisia. E-mail: lazhardhaouadi@yahoo.fr}}
\date{}
\maketitle

\begin{abstract}
In this work, we are interested by the $q$-Bessel Fourier transform with a
new approach. Many important results of this $q$-integral transform are
proved with a new constructive demonstrations and we establish in particular
the associated $q$-Fourier-Neumen expansion which involves the $q$-little
Jacobi polynomials.
\end{abstract}

\section{Introduction}

In the recent mathematical literature one finds many articles which deal
with the theory of $q$-Fourier analysis associated with the $q$-Hankel
transform. This theory was elaborated first by Koornwinder and R.F.
Swarttouw \cite{K} and then by Fitouhi and Al \cite{D1, F2}.

It should be noticed that in \cite{D1} we provided the mains results of $q$%
-Fourier analysis in particular that the $q$-Hankel transform is extended to
the $\mathcal{L}_{q,2,\nu}$ space like an isometric operator. Often we use
the crucial properties namely the positivity of the $q$-Bessel translation
operator to prove some results but these last property is not ensured for
any $q$ in the interval $]0,1[$. Thus, we will prove some main results of $q$%
-Fourier analysis without the positivity argument especially the following
statments:

- Inversion Formula in the $\mathcal{L}_{q,p,\nu}$ spaces with $p\geq 1$.

- Plancherel Formula in the $\mathcal{L}_{q,p,\nu}\cap\mathcal{L}_{q,1,\nu}$
spaces with $p>2$.

- Plancherel Formula in the $\mathcal{L}_{q,2,\nu}$ spaces.\bigskip

Note that in the paper \cite{F1} we have proved that the positivity of the $
q $-Bessel translation operator is ensured in all points of the interval $]
0,1 [$ when $\nu\geq 0$. In this article we will try to show in a clear way
the part in which the positivity of the $q$-Bessel translation operator
plays a role in $q$-Bessel Fourier analysis. In particular, when we try to prove a $
q$-version of the Young's inequality for the associated convolution.

\bigskip

Many interesting result about the uncertainty principle for the $q$-Bessel
transform was proved in the last years. We cite for examples \cite{B1,B2,B3,F3}.
There are some differences of the results cited above and our result:

In this paper the Heisenberg uncertainty inequality is established for
functions in $\mathcal{L}_{q,2,\nu}$ space.

The Hardy's inequality discuss here is a quantitative uncertainty principles
which give an information about how a function and its $q$-Bessel Fourier
transform are linked.

\bigskip

In the end of this paper we use the remarkable work in \cite{A} to establish
a new result about the $q$-Fourier-Neumen expansion involving the $q$-little
Jacobi polynomials.

\section{The $q$-Bessel transform}

The reader can see the references \cite{G,J,S} about $q$-series theory. The
references \cite{D1,F2,K} are devoted to the $q$-Bessel Fourier analysis.
Throughout this paper, we consider $0<q<1$ and $\nu >-1$. We denote by
\begin{equation*}
\mathbb{R}_{q}^{+}=\left\{ q^{n},\text{ \ }n\in \mathbb{Z}\right\} .
\end{equation*}%
The $q$-Bessel operator is defined as follows \cite{D1}
\begin{equation*}
\Delta _{q,\nu }f(x)=\frac{1}{x^{2}}\left[ f(q^{-1}x)-(1+q^{2\nu
})f(x)+q^{2\nu }f(qx)\right] .
\end{equation*}%
The eigenfunction of $\Delta _{q,\nu }$ associated with the eigenvalue $%
-\lambda ^{2}$ is the function $x\mapsto j_{\nu }(\lambda x,q^{2})$, where $%
j_{\nu }(.,q^{2})$ is the normalized $q$-Bessel function defined by \cite%
{D1,F2,G,Ko,S}
\begin{equation*}
j_{\nu }(x,q^{2})=\sum_{n=0}^{\infty }(-1)^{n}\frac{q^{n(n+1)}}{(q^{2\nu
+2},q^{2})_{n}(q^{2},q^{2})_{n}}x^{2n}.
\end{equation*}%
The $q$-Jackson integral of a function $f$ defined on $\mathbb{R^{+}}$ is
\begin{equation*}
\int_{0}^{\infty }f(t)d_{q}t=(1-q)\sum_{n\in \mathbb{Z}}q^{n}f(q^{n}).
\end{equation*}%
We denote by $\mathcal{L}_{q,p,\nu }$ the space of functions $f$ defined on $%
\mathbb{R_{q}^{+}}$ such that
\begin{equation*}
\Vert f\Vert _{q,p,\nu }=\left( \int_{0}^{\infty }|f(x)|^{p}x^{2\nu
+1}d_{q}x\right) ^{1/p}\mathrm{exist}.
\end{equation*}%
We denote by ${\mathcal{C}}_{q,0}$ the space of functions defined on $%
\mathbb{R}_{q}^{+}$ tending to $0$ as $x\rightarrow \infty $ and continuous
at $0$ equipped with the topology of uniform convergence. The space ${%
\mathcal{C}}_{q,0}$ is complete with respect to the norm
\begin{equation*}
\Vert f\Vert _{q,\infty }=\sup_{x\in \mathbb{R_{q}^{+}}}|f(x)|.
\end{equation*}%
The normalized $q$-Bessel function $j_{\nu }(.,q^{2})$ satisfies the
orthogonality relation
\begin{equation}
c_{q,\nu }^{2}\int_{0}^{\infty }j_{\nu }(xt,q^{2})j_{\nu }(yt,q^{2})t^{2\nu
+1}d_{q}t=\delta _{q}(x,y),\quad \forall x,y\in \mathbb{R}_{q}^{+}
\label{e2}
\end{equation}%
where
\begin{equation*}
\delta _{q}(x,y)=\left\{
\begin{tabular}{l}
$0$ if $x\neq y$ \\
$\frac{1}{(1-q)x^{2(\nu +1)}}$ if $x=y$%
\end{tabular}%
\ \right.
\end{equation*}%
and
\begin{equation*}
c_{q,\nu }=\frac{1}{1-q}\frac{(q^{2\nu +2},q^{2})_{\infty }}{%
(q^{2},q^{2})_{\infty }}.
\end{equation*}%
Let $f$ \ be a function defined on $\mathbb{R}_{q}^{+}$ then
\begin{equation*}
\int_{0}^{\infty }f(y)\delta _{q}(x,y)y^{2\nu +1}d_{q}y=f(x).
\end{equation*}%
The normalized $q$-Bessel function $j_{\nu }(.,q^{2})$ satisfies
\begin{equation*}
|j_{\nu }(q^{n},q^{2})|\leq \frac{(-q^{2};q^{2})_{\infty }(-q^{2\nu
+2};q^{2})_{\infty }}{(q^{2\nu +2};q^{2})_{\infty }}\left\{
\begin{array}{c}
1\quad \quad \quad \quad \quad \text{if}\quad n\geq 0 \\
q^{n^{2}-(2\nu +1)n}\quad \text{if}\quad n<0%
\end{array}%
\right. .
\end{equation*}%
The $q$-Bessel Fourier transform $\mathcal{F}_{q,\nu }$ is defined by \cite%
{D1,F2,K}
\begin{equation*}
\mathcal{F}_{q,\nu }f(x)=c_{q,\nu }\int_{0}^{\infty }f(t)j_{\nu
}(xt,q^{2})t^{2\nu +1}d_{q}t,\quad \forall x\in \mathbb{R}_{q}^{+}.
\end{equation*}

\begin{proposition}
\label{p1} Let $f\in \mathcal{L}_{q,1,\nu }$ then $\mathcal{F}_{q,\nu }f\in
\mathcal{C}_{q,0}$ and we have
\begin{equation*}
\Vert \mathcal{F}_{q,\nu }(f)\Vert _{q,\infty }\leq B_{q,\nu}\Vert f\Vert
_{q,1,v}
\end{equation*}%
where
\begin{equation*}
B_{q,\nu}=\frac{1}{1-q}\frac{(-q^{2};q^{2})_{\infty
}(-q^{2v+2};q^{2})_{\infty }}{(q^{2};q^{2})_{\infty }}.
\end{equation*}
\end{proposition}

\begin{theorem}
\label{t1} Let $f$ be a function in the $\mathcal{L}_{q,p,\nu}$ space where $%
p\geq 1$ then
\begin{equation}  \label{e1}
\mathcal{F}_{q,\nu}^{2}f=f.
\end{equation}
\end{theorem}

\begin{proof}
If $f\in \mathcal{L}_{q,p,\nu }$ then $\mathcal{F}_{q,\nu }f$ exist, and we
have
\begin{align*}
\mathcal{F}_{q,\nu }^{2}f(x)& =c_{q,\nu }\int_{0}^{\infty }\mathcal{F}%
_{q,\nu }f(t)j_{\nu }(xt,q^{2})t^{2\nu +1}d_{q}t \\
& =\int_{0}^{\infty }f(y)\left[ c_{q,\nu }^{2}\int_{0}^{\infty }j_{\nu
}(xt,q^{2})j_{\nu }(yt,q^{2})t^{2\nu +1}d_{q}t\right] y^{2\nu +1}d_{q}y \\
& =\int_{0}^{\infty }f(y)\delta _{q}(x,y)y^{2\nu +1}d_{q}y \\
& =f(x).
\end{align*}%
The computations are justified by the Fubuni's theorem: If $p>1$ then we use
the H\"{o}lder's inequality
\begin{align*}
& \int_{0}^{\infty }|f(y)|\left[ \int_{0}^{\infty }|j_{\nu }(xt,q^{2})j_{\nu
}(yt,q^{2})|t^{2\nu +1}d_{q}t\right] y^{2\nu +1}d_{q}y \\
& \leq \left[ \int_{0}^{\infty }|f(y)|^{p}y^{2\nu +1}d_{q}y\right]
^{1/p}\times \left[ \int_{0}^{\infty }\sigma (y)^{\overline{p}}y^{2\nu
+1}d_{q}y\right] ^{1/\overline{p}}.
\end{align*}%
The numbers $p$ and $\overline{p}$ above are conjugates and
\begin{equation*}
\sigma (y)=\int_{0}^{\infty }|j_{\nu }(xt,q^{2})j_{\nu }(yt,q^{2})|t^{2\nu
+1}d_{q}t,
\end{equation*}%
then
\begin{align*}
& \int_{0}^{\infty }\sigma (y)^{\overline{p}}y^{2\nu +1}d_{q}y \\
& =\int_{0}^{1}\sigma (y)^{\overline{p}}y^{2\nu +1}d_{q}y+\int_{1}^{\infty
}\sigma y^{\overline{p}}y^{2\nu +1}d_{q}y.
\end{align*}%
Note that
\begin{align*}
& \int_{0}^{1}\sigma (y)^{\overline{p}}y^{2\nu +1}d_{q}y \\
& \leq \Vert j_{\nu }(.,q^{2})\Vert _{q,\infty }^{\overline{p}}\int_{0}^{1}%
\left[ \int_{0}^{\infty }|j_{\nu }(xt,q^{2})|t^{2\nu +1}d_{q}t\right] ^{%
\overline{p}}y^{2\nu +1}d_{q}y \\
& \leq \Vert j_{\nu }(.,q^{2})\Vert _{q,\infty }^{\overline{p}}\Vert j_{\nu
}(.,q^{2})\Vert _{q,1,\nu }^{\overline{p}}x^{-2(\nu +1)\overline{p}}\left[
\int_{0}^{1}y^{2\nu +1}d_{q}y\right] <\infty ,
\end{align*}%
and
\begin{align*}
& \int_{1}^{\infty }\sigma (y)^{\overline{p}}y^{2\nu +1}d_{q}y \\
& \leq \Vert j_{\nu }(.,q^{2})\Vert _{q,\infty }^{\overline{p}}\Vert j_{\nu
}(.,q^{2})\Vert _{q,1,\nu }^{\overline{p}}\int_{1}^{\infty }\frac{y^{2\nu +1}%
}{y^{2(\nu +1)\overline{p}}}d_{q}y \\
& \leq \Vert j_{\nu }(.,q^{2})\Vert _{q,\infty }^{\overline{p}}\Vert j_{\nu
}(.,q^{2})\Vert _{q,1,\nu }^{\overline{p}}\int_{1}^{\infty }\frac{1}{%
y^{2(\nu +1)(\overline{p}-1)+1}}d_{q}y<\infty .
\end{align*}%
If $p=1$ then
\begin{align*}
& \int_{0}^{\infty }\Vert f(y)\Vert \left[ \int_{0}^{\infty }|j_{\nu
}(xt,q^{2})j_{\nu }(yt,q^{2})|t^{2\nu +1}d_{q}t\right] y^{2\nu +1}d_{q}y \\
& \leq \Vert f\Vert _{q,1,\nu }\Vert j_{\nu }(.,q^{2})\Vert _{q,\infty
}\Vert j_{\nu }(.,q^{2})\Vert _{q,1,\nu }\times \frac{1}{x^{2(\nu +1)}}.
\end{align*}
\end{proof}

\begin{theorem}
Let $f$ be a function in the $\mathcal{L}_{q,1,\nu}\cap \mathcal{L}%
_{q,p,\nu} $ space, where $p>2$ then
\begin{equation*}
\| \mathcal{F}_{q,\nu}f\| _{q,2,\nu}=\| f\|_{q,2,\nu}.
\end{equation*}
\end{theorem}

\begin{proof}
\bigskip Let $f\in \mathcal{L}_{q,1,\nu}\cap \mathcal{L}_{q,p,\nu}$ then by
Theorem \ref{t1} we see that
\begin{equation*}
\mathcal{F}_{q,\nu}^{2}f=f.
\end{equation*}
This implies
\begin{align*}
\int_{0}^{\infty }\mathcal{F}_{q,\nu}f(x)^{2}x^{2\nu+1}d_{q}x
&=\int_{0}^{\infty}\mathcal{F}_{q,\nu}f(x)\left[c_{q,\nu}\int_{0}^{\infty
}f(t)j_{\nu}(xt,q^{2})t^{2\nu+1}d_{q}t\right] x^{2\nu+1}d_{q}x \\
&=\int_{0}^{\infty }f(t)\left[c_{q,\nu}\int_{0}^{\infty }\mathcal{F }%
_{q,\nu}f(x)j_{\nu}(xt,q^{2})x^{2\nu+1}d_{q}x\right] t^{2\nu+1}d_{q}t \\
&=\int_{0}^{\infty }f(t)^{2}t^{2\nu+1}d_{q}t.
\end{align*}
The computations are justified by the Fubuni's theorem
\begin{align*}
&\int_{0}^{\infty}| f(t)|\left[c_{q,\nu}\int_{0}^{\infty}|\mathcal{F}%
_{q,\nu}f(x)||j_{\nu}(xt,q^{2})| x^{2\nu+1}d_qx\right] t^{2\nu+1}d_{q}t \\
&\leq\left[\int_{0}^{\infty }| f(t)|^{p}t^{2\nu+1}d_qt\right]^{1/p} \times %
\left[\int_{0}^{\infty}| \phi (t)| ^{\overline{p}}t^{2\nu+1}d_qt\right]^{1/%
\overline{p}},
\end{align*}
where
\begin{equation*}
\phi (t)=c_{q,\nu}\int_{0}^{\infty }| \mathcal{F}_{q,\nu}f(x) ||
j_{\nu}(xt,q^{2})| x^{2\nu+1}d_qx,
\end{equation*}
then
\begin{align*}
\| \mathcal{F}_{q,\nu}f(x)\| &\leq c_{q,\nu}\int_{0}^{\infty}| f(y)| |
j_{\nu}(xy,q^{2})|y^{2\nu+1}d_qy \\
&\leq c_{q,\nu}\left[\int_{0}^{\infty }| f(y)|^{p}y^{2\nu+1}d_qy\right]%
^{1/p}\times \left[\int_{0}^{\infty }| j_{\nu}(xy,q^{2})|^{\overline{p}%
}y^{2\nu+1}d_qy\right]^{1/\overline{p}} \\
&\leq c_{q,\nu}\left[\int_{0}^{\infty}| f(y)|^{p}y^{2\nu+1}d_{q}y\right]
^{1/p}\times\left[\int_{0}^{\infty }| j_{\nu}(y,q^{2})|^{\overline{p}%
}y^{2\nu+1}d_{q}y\right]^{1/\overline{p}}x^{-2(\nu+1)/\overline{p}} \\
&\leq c_{q,\nu}\| f\| _{q,p,\nu}\|j_{\nu}(.,q^{2})\|_{q,\overline{p}%
,\nu}x^{-2(\nu+1)/\overline{p}}.
\end{align*}
This gives
\begin{align*}
\phi (t) &\leq c_{q,\nu}^{2}\| f\| _{q,p,\nu}\|j_{\nu}(.,q^{2})\|_{q,%
\overline{p},\nu} \int_{0}^{\infty }|j_{\nu}(xt,q^{2})| x^{( 2\nu+1)
-2(\nu+1)/\overline{p}}d_qx \\
&\leq c_{q,\nu}^{2}\| f\| _{q,p,\nu}\|j_{\nu}(.,q^{2})\| _{q,\overline{p}%
,\nu}\left[\int_{0}^{\infty }| j_{\nu}(x,q^{2})| x^{2( \nu+1) /p-1}d_{q}x%
\right] t^{-2( \nu+1) /p} \\
&\leq C_{1}t^{-2( \nu+1)/p},
\end{align*}
and
\begin{align*}
\phi (t) &=c_{q,\nu}\int_{0}^{\infty }|\mathcal{F}_{q,\nu}f(x)| |
j_{\nu}(xt,q^{2})| x^{2\nu+1}d_qx \\
&=\left[c_{q,\nu}\int_{0}^{\infty }|\mathcal{F}_{q,\nu}f(x/t)| |
j_{\nu}(x,q^{2})| x^{2\nu+1}d_{q}x\right] t^{-2(\nu+1)} \\
&\leq c_{q,\nu}\| \mathcal{F}_{q,\nu}f\| _{q,\infty }\times \|
j_{\nu}(.,q^{2})\| _{q,1,\nu}\times t^{-2(\nu+1)} \\
&\leq C_{2}t^{-2(\nu+1)}.
\end{align*}
Note that
\begin{equation*}
\left\{
\begin{array}{c}
-1<-2(\nu+1)\frac{\overline{p}}{p}+2v+1 \\
-2(\nu+1)\overline{p}+2v+1<-1%
\end{array}
\right. \Leftrightarrow\left\{
\begin{array}{c}
0<-2(\nu+1)(\overline{p}-2) \\
-2(\nu+1)(\overline{p}-1)<0%
\end{array}
\right. \Leftrightarrow1<\overline{p}<2\Leftrightarrow p>2.
\end{equation*}
Hence
\begin{align*}
\int_{0}^{\infty }| \phi (t)| ^{\overline{p}}t^{2\nu+1}d_{q}t
&=\int_{0}^{1}| \phi (t)| ^{\overline{p} }t^{2\nu+1}d_{q}t+\int_{1}^{\infty
}| \phi (t)| ^{\overline{p}}t^{2\nu+1}d_{q}t \\
&\leq C_{1}\int_{0}^{1}t^{-2( v+1) \overline{p}/p}t^{2\nu+1}d_{q}t+C_{2}%
\int_{1}^{\infty }t^{-2( v+1) \overline{p}}t^{2\nu+1}d_{q}t<\infty,
\end{align*}
which prove the result.
\end{proof}

\begin{theorem}
\label{t3} Let $f$ be a function in the $\mathcal{L}_{q,2,\nu}$ space then
\begin{equation*}
\| \mathcal{F}_{q,\nu}f\|_{q,2,\nu}=\| f\|_{q,2,\nu}.
\end{equation*}
\end{theorem}

\begin{proof}
We introduce the function $\psi _{x}$ as follows
\begin{equation*}
\psi _{x}(t)=c_{q,\nu }j_{\nu }(tx,q^{2}).
\end{equation*}%
The inner product $\langle ,\rangle $ in the Hilbert space $\mathcal{L}%
_{q,2,\nu }$ is defined by
\begin{equation}
f,g\in \mathcal{L}_{q,2,\nu }\Rightarrow \langle f,g\rangle
=\int_{0}^{\infty }f(t)g(t)t^{2\nu +1}d_{q}t.  \label{e4}
\end{equation}%
Using (\ref{e2}) we write
\begin{equation*}
x\neq y\Rightarrow \langle \psi _{x},\psi _{y}\rangle =0
\end{equation*}%
\begin{equation*}
\Vert \psi _{x}\Vert _{q,2,\nu }^{2}=\frac{1}{1-q}x^{-2(\nu +1)}.
\end{equation*}%
We have
\begin{equation*}
\mathcal{F}_{q,\nu }f(x)=\langle f,\psi _{x}\rangle ,
\end{equation*}%
and by Theorem \ref{t1}
\begin{equation*}
f\in \mathcal{L}_{q,2,\nu }\Rightarrow \mathcal{F}_{q,\nu }^{2}f=f,
\end{equation*}%
then
\begin{equation*}
\langle f,\psi _{x}\rangle =0,\forall x\in \mathbb{R}_{q}^{+}\mathbb{%
\Rightarrow }\mathcal{F}_{q,\nu }f(x)=0,\forall x\in \mathbb{R}_{q}^{+}%
\mathbb{\Rightarrow }f=0.
\end{equation*}%
Hence, $\{\psi _{x},x\in \mathbb{R}_{q}^{+}\mathbb{\}}$ form an orthogonal
basis of the Hilbert space $\mathcal{L}_{q,2,\nu }$ and we have
\begin{equation*}
\overline{\{\psi _{x},\text{ \ \ }\forall x\in \mathbb{R}_{q}^{+}\mathbb{\}}}%
=\mathcal{L}_{q,2,\nu }.
\end{equation*}%
Now
\begin{equation*}
f\in \mathcal{L}_{q,2,\nu }\Rightarrow f=\sum_{x\in \mathbb{R_{q}^{+}}}\frac{%
1}{\Vert \psi _{x}\Vert _{q,2,\nu }^{2}}\langle f,\psi _{x}\rangle \psi _{x},
\end{equation*}%
and then
\begin{equation*}
\Vert f\Vert _{q,2,\nu }^{2}=\sum_{x\in \mathbb{R_{q}^{+}}}\frac{1}{\Vert
\psi _{x}\Vert _{q,2,\nu }^{2}}\langle f,\psi _{x}\rangle
^{2}=(1-q)\sum_{x\in \mathbb{R_{q}^{+}}}x^{2(\nu +1)}\mathcal{F}_{q,\nu
}f(x)^{2}=\Vert \mathcal{F}_{q,\nu }f\Vert _{q,2,\nu }^{2},
\end{equation*}%
which achieve the proof.
\end{proof}

\begin{proposition}
\label{p2} Let $f\in \mathcal{L}_{q,p,\nu }$ where $p\geq 1$ then $\mathcal{F}%
_{q,\nu }f\in \mathcal{L}_{q,\overline{p},\nu }$. If  $1\leq p\leq 2$
then
\begin{equation}\label{e12}
\Vert \mathcal{F}_{q,\nu }f\Vert _{q,\overline{p},\nu }\leq B_{q,\nu }^{%
\frac{2}{p}-1}\Vert f\Vert _{q,p,\nu }.
\end{equation}
\end{proposition}

\begin{proof}
This is an immediate consequence of Proposition \ref{p1}, Theorem \ref{t3},
the Riesz-Thorin theorem and the inversion formula (\ref{e1}).
\end{proof}

The $q$-translation operator is given as follow
\begin{equation*}
T_{q,x}^{\nu }f(y)=c_{q,\nu }\int_{0}^{\infty }\mathcal{F}_{q,\nu
}f(t)j_{\nu }(yt,q^{2})j_{\nu }(xt,q^{2})t^{2\nu +1}d_{q}t\bigskip .
\end{equation*}%
Let us now introduce
\begin{equation*}
Q_{\nu}=\left\{ q\in ]0,1[,\quad T_{q,x}^{\nu}\quad \text{is positive for all%
}\quad x\in \mathbb{R}_{q}^{+}\right\}
\end{equation*}%
the set of the positivity of $T_{q,x}^{\nu}$. We recall that $T_{q,x}^{\nu}$
is called positive if $T_{q,x}^{v}f\geq 0$ for $f\geq 0$. In a recent paper
\cite{D2} it was proved that if $-1<\nu<\nu^{\prime }$ then $Q_{\nu}\subset
Q_{\nu^{\prime }}$ . As a consequence :

\begin{description}
\item[-] If $0\leq \nu$ then $Q_{\nu}=]0,1[.$

\item[-] If $-\frac{1}{2}\leq \nu<0$ then $]0,q_{0}]\subset Q_{-\frac{1}{2}%
}\subset $ $Q_{\nu}\subsetneq]0,1[$,\quad $q_0\simeq 0.43$.

\item[-] If $-1<\nu\leq -\frac{1}{2}$ then $Q_{\nu}\subset Q_{-\frac{1}{2}}.$
\end{description}

\begin{theorem}
Let $f\in \mathcal{L}_{q,p,\nu }$ then $T_{q,x}^{\nu }f$ exists and we have
\begin{equation*}
\int_{0}^{\infty }T_{q,x}^{\nu }f(y)y^{2\nu +1}d_{q}y=\int_{0}^{\infty
}f(y)y^{2\nu +1}d_{q}y.
\end{equation*}%
and
\begin{equation*}
T_{q,x}^{\nu }f(y)=\int_{0}^{\infty }f(z)D_{\nu }(x,y,z)z^{2\nu +1}d_{q}z,
\end{equation*}%
where
\begin{equation*}
D_{\nu }(x,y,z)=c_{q,\nu }^{2}\int_{0}^{\infty }j_{\nu }(xs,q^{2}j_{\nu
}(ys,q^{2}j_{\nu }(zs,q^{2})s^{2\nu +1}d_{q}s.
\end{equation*}%
If we suppose that $T_{q,x}^{\nu }$ is a positive operator then for all $p\geq 1$ we
have
\begin{equation}\label{e9}
\Vert T_{q,x}^{\nu }f\Vert _{q,p,\nu }\leq \Vert f\Vert _{q,p,\nu }.
\end{equation}
\end{theorem}

\begin{proof}
We write the operator $T_{q,x}^{\nu}$ in the following form
\begin{eqnarray*}
T_{q,x}^{\nu}f(y) &=&c_{q,\nu}\int_{0}^{\infty }\mathcal{F}%
_{q,v}f(z)j_{v}(xz,q^{2})j_{v}(yz,q^{2})z^{2v+1}d_{q}z \\
&=&\mathcal{F}_{q,\nu}\left[ \mathcal{F}_{q,\nu}f(z)j_{\nu}(xz,q^{2})\right]
(y).
\end{eqnarray*}%
So we have
\begin{eqnarray*}
\int_{0}^{\infty }T_{q,x}^{\nu}f(y)y^{2\nu+1}d_{q}y &=&\int_{0}^{\infty }%
\mathcal{F}_{q,\nu}\left[ \mathcal{F}_{q,\nu}f(z)j_{\nu}(xz,q^{2})\right]
(y)y^{2\nu+1}d_{q}y \\
&=&\frac{1}{c_{q,\nu}}c_{q,\nu}\int_{0}^{\infty }\mathcal{F}_{q,\nu}\left[
\mathcal{F}_{q,\nu}f(z)j_{\nu}(xz,q^{2})\right] (y)j_{\nu}(0,q^{2})y^{2%
\nu+1}d_{q}y \\
&=&\frac{1}{c_{q,\nu}}\mathcal{F}_{q,\nu}^{2}\left[ \mathcal{F}%
_{q,\nu}f(z)j_{\nu}(xz,q^{2})\right] (0) \\
&=&\frac{1}{c_{q,\nu}}\mathcal{F}_{q,\nu}f(0) \\
&=&\int_{0}^{\infty }f(y)y^{2v+1}d_{q}y.
\end{eqnarray*}%
On the other hand
\begin{eqnarray*}
T_{q,x}^{\nu}f(y) &=&c_{q,\nu}\int_{0}^{\infty }\mathcal{F}%
_{q,\nu}f(z)j_{\nu}(xz,q^{2})j_{\nu}(yz,q^{2})z^{2\nu+1}d_{q}z \\
&=&c_{q,\nu}\int_{0}^{\infty }\left[ c_{q,\nu}\int_{0}^{\infty
}f(t)j_{\nu}(tz,q^{2})t^{2\nu+1}d_{q}t\right] j_{\nu}(xz,q^{2})j_{%
\nu}(yz,q^{2})z^{2v+1}d_{q}z \\
&=&\int_{0}^{\infty }\left[ c_{q,\nu}^{2}\int_{0}^{\infty
}j_{\nu}(xz,q^{2})j_{\nu}(yz,q^{2})j_{\nu}(tz,q^{2})z^{2\nu+1}d_{q}z\right]
f(t)t^{2\nu+1}d_{q}t \\
&=&\int_{0}^{\infty }D_{q,\nu }(x,y,t)f(t)t^{2\nu+1}d_{q}t.
\end{eqnarray*}%
The computations are justified by the Fubuni's theorem
\begin{eqnarray*}
&&\int_{0}^{\infty }\left[ \int_{0}^{\infty }\left\vert f(t)\right\vert
\left\vert j_{\nu}(tz,q^{2})\right\vert t^{2\nu+1}d_{q}t\right] \left\vert
j_{v}(xz,q^{2})\right\vert \left\vert j_{\nu}(yz,q^{2})\right\vert
z^{2v+1}d_{q}z \\
&\leq &\left\Vert f\right\Vert _{q,p,\nu}\int_{0}^{\infty }\left[
\int_{0}^{\infty }\left\vert j_{\nu}(tz,q^{2})\right\vert ^{\overline{p}%
}t^{2\nu+1}d_{q}t\right] ^{\frac{1}{\overline{p}}}\left\vert
j_{\nu}(xz,q^{2})\right\vert \left\vert j_{v}(yz,q^{2})\right\vert
z^{2\nu+1}d_{q}z \\
&\leq &\left\Vert f\right\Vert _{q,p,\nu}\left\Vert
j_{\nu}(.,q^{2})\right\Vert _{q,\overline{p},\nu}\int_{0}^{\infty
}\left\vert j_{\nu}(xz,q^{2})\right\vert \left\vert
j_{\nu}(yz,q^{2})\right\vert z^{2(\nu+1)\left( 1-\frac{1}{\overline{p}}%
\right) -1}d_{q}z.
\end{eqnarray*}%
Now suppose that $T_{q,x}^{\nu }$ is positive. Given a function $f\in
\mathcal{C}_{q,0}$ we obtains
\begin{eqnarray*}
\left\vert T_{q,x}^{\nu}f(y)\right\vert &=&\left\vert \int_{0}^{\infty
}D_{q,\nu}(x,y,t)f(t)t^{2\nu+1}d_{q}t\right\vert \\
&\leq &\int_{0}^{\infty }\left\vert D_{q,\nu}(x,y,t)\right\vert \left\vert
f(t)\right\vert t^{2v+1}d_{q}t \\
&\leq &\left[ \int_{0}^{\infty }D_{q,\nu}(x,y,t)t^{2v+1}d_{q}t\right]
\left\Vert f\right\Vert _{q,\infty }=\left\Vert f\right\Vert _{q,\infty }
\end{eqnarray*}%
which implies
\begin{equation*}
\left\Vert T_{q,x}^{\nu}f\right\Vert _{q,\infty }\leq \left\Vert
f\right\Vert _{q,\infty }.
\end{equation*}%
If the function $f\in \mathcal{L}_{q,1,\nu }$ then we obtains
\begin{eqnarray*}
\left\Vert T_{q,x}^{\nu}f\right\Vert _{q,1,\nu} &=&\int_{0}^{\infty
}\left\vert T_{q,x}^{\nu}f(y)\right\vert y^{2\nu+1}d_{q}y \\
&\leq &\int_{0}^{\infty }\left[ \int_{0}^{\infty }\left\vert
D_{q,\nu}(x,y,t)\right\vert \left\vert f(t)\right\vert t^{2\nu+1}d_{q}t%
\right] y^{2\nu+1}d_{q}y \\
&\leq &\int_{0}^{\infty }\left[ \int_{0}^{\infty
}D_{q,\nu}(x,y,t)y^{2\nu+1}d_{q}y\right] \left\vert f(t)\right\vert
t^{2\nu+1}d_{q}t \\
&\leq &\int_{0}^{\infty }\left\vert f(t)\right\vert
t^{2\nu+1}d_{q}t=\left\Vert f\right\Vert _{q,1,\nu}.
\end{eqnarray*}
The result is a consequence of the Riesz-Thorin theorem.

\bigskip

Notice that the kernel $D_{q,\nu }(x,y,t)$ can be written as follows
\begin{eqnarray*}
D_{q,\nu }(x,y,t) &=&c_{q,\nu }^{2}\int_{0}^{\infty }j_{\nu
}(xz,q^{2})j_{\nu }(yz,q^{2})j_{\nu }(tz,q^{2})z^{2\nu +1}d_{q}z \\
&=&c_{q,\nu }\mathcal{F}_{q,\nu }\left[ j_{\nu }(xz,q^{2})j_{\nu }(yz,q^{2})%
\right] (t),
\end{eqnarray*}%
which implies
\begin{eqnarray*}
\int_{0}^{\infty }D_{q,\nu }(x,y,t)t^{2\nu +1}d_{q}t &=&c_{q,\nu
}\int_{0}^{\infty }\mathcal{F}_{q,\nu }\left[ j_{\nu }(xz,q^{2})j_{\nu
}(yz,q^{2})\right] (t)t^{2\nu +1}d_{q}t \\
&=&\mathcal{F}_{q,\nu }^{2}\left[ j_{\nu }(xz,q^{2})j_{\nu }(yz,q^{2})\right]
(0)=1.
\end{eqnarray*}
\end{proof}

\bigskip

The $q$-convolution product is defined by
\begin{equation*}
f\ast _{q}g=\mathcal{F}_{q,\nu }\left[ \mathcal{F}_{q,\nu }f\times \mathcal{F%
}_{q,\nu }g\right].
\end{equation*}

\begin{theorem}
Let $\ 1\leq p,r,s$ such that%
\begin{equation*}
\frac{1}{p}+\frac{1}{r}-1=\frac{1}{s}
\end{equation*}%
Given two functions $f\in \mathcal{L}_{q,p,\nu }$ and $g\in \mathcal{L}%
_{q,r,\nu }$ then $f\ast _{q}g$ exists and we have
\begin{equation*}
f\ast _{q}g(x)=c_{q,\nu }\int_{0}^{\infty }T_{q,x}^{\nu }f(y)g(y)y^{2\nu
+1}d_{q}y.
\end{equation*}%
and%
\begin{equation*}
f\ast _{q}g \in \mathcal{L}_{q,s,\nu }.
\end{equation*}%
\begin{equation*}
\mathcal{F}_{q,\nu }(f\ast _{q}g)=\mathcal{F}_{q,\nu }(f)\times \mathcal{F}%
_{q,\nu }(g).
\end{equation*}
If $s\geq 2$ then
\begin{equation}\label{e10}
\left\Vert f\ast _{q}g\right\Vert _{q,s,\nu }\leq B_{q,\nu }\left\Vert
f\right\Vert _{q,p,\nu }\left\Vert g\right\Vert _{q,r,\nu }.
\end{equation}
If we suppose that $T_{q,x}^{\nu }$ is a positive operator then
\begin{equation}\label{e11}
\left\Vert f\ast _{q}g\right\Vert _{q,s,\nu }\leq c_{q,\nu }\left\Vert
f\right\Vert _{q,p,\nu }\left\Vert g\right\Vert _{q,r,\nu }.
\end{equation}
\end{theorem}

\begin{proof}
We have
\begin{eqnarray*}
f\ast _{q}g(x) &=&\mathcal{F}_{q,\nu }\left[ \mathcal{F}_{q,\nu }f\times
\mathcal{F}_{q,\nu }g\right] (x) \\
&=&c_{q,\nu }\int_{0}^{\infty }\mathcal{F}_{q,\nu }f(y)\times \mathcal{F}%
_{q,\nu }g(y)j_{\nu }(xy,q^{2})y^{2\nu +1}d_{q}y \\
&=&c_{q,\nu }\int_{0}^{\infty }\mathcal{F}_{q,\nu }f(y)\times \left[
c_{q,\nu }\int_{0}^{\infty }g(z)j_{\nu }(zy,q^{2})z^{2\nu +1}d_{q}z\right]
j_{v}(xy,q^{2})y^{2\nu +1}d_{q}y \\
&=&c_{q,\nu }\int_{0}^{\infty }\left[ c_{q,\nu }\int_{0}^{\infty }\mathcal{F}%
_{q,\nu }f(y)j_{v}(zy,q^{2})j_{\nu }(xy,q^{2})y^{2\nu +1}d_{q}y\right]
g(z)z^{2\nu +1}d_{q}z \\
&=&c_{q,\v}\int_{0}^{\infty }T_{q,x}^{\nu }f(z)g(z)z^{2\nu +1}d_{q}z.
\end{eqnarray*}%
The computations are justified by the Fubuni's theorem
\begin{eqnarray*}
&&\int_{0}^{\infty }\left\vert F_{q,\nu }f(y)\right\vert \times \left[
\int_{0}^{\infty }\left\vert g(z)\right\vert \times \left\vert j_{\nu
}(zy,q^{2})\right\vert z^{2\nu +1}d_{q}z\right] \left\vert j_{\nu
}(xy,q^{2})\right\vert y^{2\nu +1}d_{q}y \\
&\leq &\left\Vert g\right\Vert _{q,r,\nu }\int_{0}^{\infty }\left\vert
F_{q,\nu }f(y)\right\vert \times \left[ \int_{0}^{\infty }\left\vert j_{\nu
}(zy,q^{2})\right\vert ^{\overline{r}}z^{2\nu +1}d_{q}z\right] ^{\frac{1}{%
\overline{r}}}\left\vert j_{\nu }(xy,q^{2})\right\vert y^{2\nu +1}d_{q}y \\
&\leq &\left\Vert g\right\Vert _{q,r,\nu }\left\Vert j_{\nu
}(.,q^{2})\right\Vert _{q,\overline{r},\nu }\int_{0}^{\infty }\left\vert
F_{q,\nu }f(y)\right\vert \times \left[ \left\vert
j_{v}(xy,q^{2})\right\vert y^{-\frac{2\nu +2}{\overline{r}}}\right] y^{2\nu
+1}d_{q}y \\
&\leq &\left\Vert g\right\Vert _{q,r,\nu }\left\Vert j_{\nu
}(.,q^{2})\right\Vert _{q,\overline{r},\nu }\left\Vert F_{q,\nu
}f\right\Vert _{q,\overline{p},v}\left( \int_{0}^{\infty }\left[ \left\vert
j_{\nu }(xy,q^{2})\right\vert y^{-\frac{2\nu +2}{\overline{r}}}\right]
^{p}y^{2\nu +1}d_{q}y\right) ^{\frac{1}{p}} \\
&\leq &\left\Vert g\right\Vert _{q,r,\nu }\left\Vert j_{\nu
}(.,q^{2})\right\Vert _{q,\overline{r},\nu }\left\Vert F_{q,\nu
}f\right\Vert _{q,\overline{p},\nu }\left( \int_{0}^{\infty }\left\vert
j_{\nu }(xy,q^{2})\right\vert ^{p}y^{2(\nu +1)\left( 1-\frac{p}{\overline{r}}%
\right) -1}d_{q}y\right) ^{\frac{1}{p}}.
\end{eqnarray*}%
From Proposition \ref{p2} we deduce that
\begin{equation*}
\mathcal{F}_{q,\nu }f\in \mathcal{L}_{q,\overline{p},\nu }\text{ }\mathrm{and%
}\text{ }\mathcal{F}_{q,\nu }g\in \mathcal{L}_{q,\overline{r},\nu }.
\end{equation*}%
Then, using the H\"{o}lder inequality and the fact that%
\begin{equation*}
\frac{1}{\overline{p}}+\frac{1}{\overline{r}}=\frac{1}{\overline{s}}
\end{equation*}%
to conclude that
\begin{equation*}
\mathcal{F}_{q,\nu }f\times \mathcal{F}_{q,\nu }g\in \mathcal{L}_{q,%
\overline{s},\nu }.
\end{equation*}%
Which implies that
\begin{equation*}
f\ast _{q}g=\mathcal{F}_{q,\nu }\left[ \mathcal{F}_{q,\nu }f\times \mathcal{F%
}_{q,\nu }g\right] \in \mathcal{L}_{q,s,\nu }
\end{equation*}%
and by the inversion formula (\ref{e1}) we obtain
\begin{equation*}
\mathcal{F}_{q,\nu }\left( f\ast _{q}g\right) =\mathcal{F}_{q,\nu }f\times
\mathcal{F}_{q,\nu }g.
\end{equation*}%
Suppose that $s\geq 2,$ so $1\leq \overline{s}\leq 2$ and we can write%
\begin{eqnarray*}
\left\Vert f\ast _{q}g\right\Vert _{q,s,\nu } &=&\left\Vert \mathcal{F}%
_{q,\nu }\left[ \mathcal{F}_{q,\nu }f\times \mathcal{F}_{q,\nu }g\right]
\right\Vert _{q,s,\nu } \\
&\leq &B_{q,\nu }^{\frac{2}{\overline{s}}-1}\left\Vert \mathcal{F}_{q,\nu
}f\right\Vert _{q,\overline{p},\nu }\left\Vert \mathcal{F}_{q,\nu
}g\right\Vert _{q,\overline{r},\nu } \\
&\leq &B_{q,\nu }^{\frac{2}{\overline{s}}-1}B_{q,\nu }^{\frac{2}{p}%
-1}B_{q,\nu }^{\frac{2}{r}-1}\left\Vert f\right\Vert _{q,p,\nu }\left\Vert
g\right\Vert _{q,r,\nu } \\
&\leq &B_{q,\nu }\left\Vert f\right\Vert _{q,p,\nu }\left\Vert g\right\Vert
_{q,r,\nu }.
\end{eqnarray*}%
Now suppose that $T_{q,x}^{\nu }$ is a positive operator.

We introduce the operator $K_{f}$ as follows%
\begin{equation*}
K_{f}g(x)=c_{q,\v}\int_{0}^{\infty }T_{q,x}^{\nu }f(z)g(z)z^{2\nu +1}d_{q}z.
\end{equation*}
By the H\"{o}lder inequality and (\ref{e9}) we get
\begin{equation*}
\left\Vert K_{f}g\right\Vert _{q,\infty }\leq c_{q,\nu }\left\Vert
f\right\Vert _{q,p,\nu }\left\Vert g\right\Vert _{q,\overline{p},\nu }.
\end{equation*}
The Minkowski inequality leads to
\begin{equation*}
\left\Vert K_{f}g\right\Vert _{q,p,\nu }\leq c_{q,\nu }\left\Vert
f\right\Vert _{q,p,\nu }\left\Vert g\right\Vert _{q,1,\nu }.
\end{equation*}
Hence we have
\begin{equation*}
K_{f}:\mathcal{L}_{q,\overline{p},\nu }\rightarrow \mathcal{C}_{q,0},
\quad K_{f}:\mathcal{L}_{q,1,\nu }\rightarrow \mathcal{L}_{q,p,\nu }.
\end{equation*}
Then the operator $K_f$ satisfies
\begin{equation*}
K_{f}:\mathcal{L}_{q,r,\nu }\rightarrow \mathcal{L}_{q,s,\nu }
\end{equation*}
and
\begin{equation*}
\left\Vert f\ast _{q}g\right\Vert _{q,s,\nu }=\|K_{f}g\|_{q,s,\nu }\leq c_{q,\nu }\left\Vert
f\right\Vert _{q,p,\nu }\left\Vert g\right\Vert _{q,r,\nu }.
\end{equation*}
\end{proof}

\begin{remark}
We discuss  here the  sharp  results  for the  Hausdorf-Young  inequality provided above.
An inequality already sharper than (\ref{e10}) is given in formula (\ref{e11}).
In fact we have $c_{q,\nu}<B_{q,\nu}$.

To  obtained (\ref{e11}) without the positivity argument, we can do  by using which is a $q$-Riemann-Liouville
fractional integral generalizing the $q$-Mehler integral representation for the $q$-Bessel  function $j_\nu(.,q^2)$
which can be proved in a straightforward way  \cite{F2}
$$
j_{\v}(\lambda,q^{2})=[2\v]_q\int_0^1 \frac{(q^2t^2,q^{2})_\I}{(q^{2\v}t^2,q^{2})_\I}j_0(\lambda t,q^{2})td_qt
$$
together with the inequalities for the $q$-Bessel  function which
is given as formula (24) in the paper \cite{B3}
$$
|j_0 (x;q^2)| \leq 1,\quad\forall  x\in {\mathbb R}^+_q. $$
Combine this formulas we arrive at
$$
|j_\nu (x;q^2)| \leq 1,\quad\forall  x\in {\mathbb R}^+_q,\quad \v\geq 0.
$$
Then the inequalities (\ref{e12}) can be written as follows
$$
\Vert \mathcal{F}_{q,\nu }f\Vert _{q,\overline{p},\nu }\leq c_{q,\nu }^{
\frac{2}{p}-1}\Vert f\Vert _{q,p,\nu }.
$$
This should give the sharpest version of (\ref{e10}) in the cases $\v\geq 0$. Unfortunately
the positivity of the operator $T_{q,x}^\v$ is satisfied in this case.

In fact we can prove that if we are in the positivity cases then
\begin{equation*}
\left\Vert j_{\v}(.,q^{2})\right\Vert _{q,\infty }\leq 1.
\end{equation*}
To prove this recalling that
\begin{equation*}
T_{q,x}^{\v}j_{\v}(y,q^{2})=j_{\v}(x,q^{2})j_{\v}(y,q^{2}).
\end{equation*}
So we have
\begin{equation*}
\int_{0}^{\infty
}D_{\v}(x,y,t)j_{v}(t,q^{2})t^{2v+1}d_{q}t=j_{\v}(x,q^{2})j_{\v}(y,q^{2}).
\end{equation*}
We obtains for all $x,y\in{\mathbb R}_q^+$
\begin{eqnarray*}
\left\vert j_{\v}(x,q^{2})\right\vert \times \left\vert
j_{\v}(y,q^{2})\right\vert &\leq &\int_{0}^{\infty }D_{\v}(x,y,t)\left\vert
j_{\v}(t,q^{2})\right\vert t^{2\v+1}d_{q}t \\
&\leq &\left[ \int_{0}^{\infty }D_{v}(x,y,t)t^{2\v+1}d_{q}t\right] \left\Vert
j_{\v}(.,q^{2})\right\Vert _{q,\infty }.
\end{eqnarray*}
The fact that
$$
\int_{0}^{\infty }D_{\v}(x,y,t) t^{2\v+1}d_{q}t=1
$$
implies
\begin{equation*}
\left\Vert j_{\v}(.,q^{2})\right\Vert _{q,\infty }^{2}\leq \left\Vert
j_{\v}(.,q^{2})\right\Vert _{q,\infty }
\end{equation*}
which gives the result.

\end{remark}

\section{Uncertainty principle}

We introduce two $q$-difference operators
\begin{equation*}
\partial _{q}f(x)=\frac{f(q^{-1}x)-f(x)}{x}
\end{equation*}
and
\begin{equation*}
\partial _{q}^{\ast }f(x)=\frac{f(x)-q^{2\nu+1}f(qx)}{x}.
\end{equation*}
Then we have
\begin{equation*}
\partial _{q}\partial _{q}^{\ast }f(x)=\partial _{q}^{\ast }\partial
_{q}f(x)=\Delta _{q,\nu}f(x).
\end{equation*}

\begin{proposition}
If $\left\langle \partial _{q}f,g\right\rangle $ exist \ and $\underset{%
a\rightarrow \infty }{\lim }\left\vert a^{2\nu+1}f(q^{-1}a)g(a)\right\vert
=0 $ then
\begin{equation*}
\left\langle \partial _{q}f,g\right\rangle =-\left\langle f,\partial
_{q}^{\ast }g\right\rangle .
\end{equation*}
\end{proposition}

\begin{proof}
The following computation
\begin{eqnarray*}
&&\int_{0}^{a}\partial _{q}f(x)g(x)x^{2\nu+1}d_{q}x \\
&=&\int_{0}^{a}\frac{f(q^{-1}x)-f(x)}{x}g(x)x^{2\nu+1}d_{q}x \\
&=&\int_{0}^{a}\frac{f(q^{-1}x)}{x}g(x)x^{2\nu+1}d_{q}x-\int_{0}^{a}\frac{%
f(x)}{x}g(x)x^{2\nu+1}d_{q}x \\
&=&q^{2\nu+1}\int_{0}^{q^{-1}a}\frac{f(x)}{x}g(qx)x^{2\nu+1}d_{q}x-%
\int_{0}^{a}\frac{f(x)}{x}\partial _{q}g(x)x^{2\nu+1}d_{q}x \\
&=&q^{2\nu+1}\int_{0}^{a}\frac{f(x)}{x}\partial
_{q}g(qx)x^{2\nu+1}d_{q}x-\int_{0}^{a}\frac{f(x)}{x}g(x)x^{2v+1}d_{q}x+a^{2%
\nu+1}f(q^{-1}a)g(a) \\
&=&-\int_{0}^{a}f(x)\frac{g(x)-q^{2\nu+1}g(qx)}{x}x^{2\nu+1}d_{q}x+a^{2%
\nu+1}f(q^{-1}a)g(a) \\
&=&=-\int_{0}^{a}f(x)\partial _{q}^{\ast
}g(x)x^{2\nu+1}d_{q}x+a^{2\nu+1}f(q^{-1}a)g(a)
\end{eqnarray*}%
leads to the result.
\end{proof}

\begin{corollary}\label{c2}
If $f\in \mathcal{L}_{q,2,\nu}$ such that $x\mathcal{F}_{q,\nu}f\in \mathcal{%
L} _{q,2,\nu}$ then
\begin{equation*}
\left\Vert \partial _{q}f\right\Vert _{2}=\left\Vert x\mathcal{F}%
_{q,\nu}f\right\Vert _{2}.
\end{equation*}
\end{corollary}

\begin{proof}
In fact we have
\begin{eqnarray*}
\left\Vert \partial _{q}f\right\Vert _{2}^{2} &=&\left\langle \partial
_{q}f,\partial _{q}f\right\rangle =-\left\langle f,\partial _{q}^{\ast
}\partial _{q}f\right\rangle \\
&=&-\left\langle f,\Delta _{q,\nu}f\right\rangle \\
&=&-\left\langle \mathcal{F}_{q,\nu}f,\mathcal{F}_{q,\nu}\Delta
_{q,\nu}f\right\rangle \\
&=&\left\langle \mathcal{F}_{q,\nu}f,x^{2}\mathcal{F}_{q,\nu}f\right\rangle
\\
&=&\left\Vert x\mathcal{F}_{q,\nu}f\right\Vert _{2}^{2},
\end{eqnarray*}
which prove the result.
\end{proof}

\begin{theorem}
Assume that $f$  belongs to the space $\mathcal{L}_{q,2,\nu}.$
Then the $q$-Bessel transform satisfies the
following uncertainty principal
\begin{equation*}
\left\Vert f\right\Vert _{2}^{2}\leq k_{q,v}\left\Vert xf\right\Vert
_{2}\left\Vert x\mathcal{F}_{q,\nu}f\right\Vert _{2}
\end{equation*}
where
\begin{equation*}
k_{q,\nu}=\frac{\left[1+\sqrt{q}\times q^{\nu+1}\right]}{1-q^{2(\nu+1)}}.
\end{equation*}
\end{theorem}

\begin{proof}
In fact
\begin{equation*}
\partial _{q}^{\ast }xf=f(x)-q^{2\nu+2}f(qx)
\end{equation*}
\begin{equation*}
x\partial _{q}f=f(q^{-1}x)-f(x).
\end{equation*}
We introduce the following operator
\begin{equation*}
\Lambda _{q}f(x)=f(qx),
\end{equation*}
then
\begin{equation*}
\left\langle \Lambda _{q}f,g\right\rangle =q^{-2(\nu+1)}\left\langle
f,\Lambda _{q}^{-1}g\right\rangle .
\end{equation*}
So
\begin{equation*}
\frac{1}{1-q^{2(\nu+1)}}\left[ \partial _{q}^{\ast }xf(x)-q^{2\nu+2}\Lambda
_{q}x\partial _{q}f(x)\right] =f(x)
\end{equation*}
Assume that $xf$ and $x\mathcal{F}_{q,\nu}f$ belongs to the space $\mathcal{L}_{q,2,\nu}$. Then  we have
\begin{equation*}
\left\langle f,f\right\rangle =-\frac{1}{1-q^{2(\nu+1)}}\left\langle
xf,\partial _{q}f\right\rangle -\frac{1}{1-q^{2(\nu+1)}}\left\langle
\partial _{q}f,x\Lambda _{q}^{-1}f\right\rangle .
\end{equation*}
By Cauchy-Schwartz inequality we get
\begin{equation*}
\left\langle f,f\right\rangle \leq \frac{1}{1-q^{2(\nu+1)}} \left\Vert
xf\right\Vert _{2}\left\Vert \partial _{q}f\right\Vert _{2}+\frac{1}{%
1-q^{2(\nu+1)}}\left\Vert \partial _{q}f\right\Vert _{2}\left\Vert x\Lambda
_{q}^{-1}f\right\Vert _{2}.
\end{equation*}
On the other hand
\begin{equation*}
\left\Vert x\Lambda _{q}^{-1}f\right\Vert _{2}=\sqrt{q}\times
q^{\nu+1}\left\Vert xf\right\Vert _{2},
\end{equation*}
Corollary \ref{c2} leads to the result.
\end{proof}

\section{Hardy's theorem}

The following Lemma from complex analysis is crucial for the proof of our
main theorem.

\begin{lemma}
\label{l1} For every $p\in\mathbb{N}$, there exist $\sigma_p>0$ for which
\begin{equation*}
|z|^{2p}|j_{\nu}(z,q^2)|<\sigma_p e^{|z|},\quad\forall z\in\mathbb{C}.
\end{equation*}
\end{lemma}

\begin{proof}
In fact
\begin{align*}
|z|^{2p}|j_{\nu}(z,q^2)|&\leq\frac{1}{(q^2,q^2)_{\infty}
(q^{2\nu+2},q^2)_{\infty}}\sum_{n=0}^\infty q^{n(n-1)}|z|^{2n+2p} \\
&\leq\frac{q^{p(p+1)}}{(q^2,q^2)_{\infty} (q^{2\nu+2},q^2)_{\infty}}%
\sum_{n=p}^\infty q^{n(n-2p-1)}|z|^{2n}.
\end{align*}
Now using the Stirling's formula
\begin{equation*}
n!\sim\sqrt{2\pi n}\frac{n^n}{e^n},
\end{equation*}
we see that there exist an entire $n_0\geq p$ such that
\begin{equation*}
q^{n(n-2p-1)}<\frac{1}{(2n)!},\quad\forall n\geq n_0,
\end{equation*}
which implies
\begin{equation*}
\sum_{n=n_0}^\infty q^{n(n-2p-1)}|z|^{2n}<\sum_{n=n_0}^\infty \frac{1}{(2n)!}
|z|^{2n}<e^{|z|}.
\end{equation*}
Finally there exist $\sigma_p>0$ such that
\begin{equation*}
\frac{|z|^{2p}|j_{\nu}(z,q^2)|}{e^{|z|}}<\sigma_p,\quad\forall z\in\mathbb{C}
\end{equation*}
This complete the proof.
\end{proof}

\begin{lemma}
\label{l2} Let $h$ be an entire function on $\mathbb{C}$ such that
\begin{equation*}
|h(z)|\leq C e^{a|z|^2},\quad z\in\mathbb{C},
\end{equation*}

\begin{equation*}
|h(x)|\leq Ce^{-ax^{2}},\quad x\in \mathbb{R},
\end{equation*}
for some positive constants $a$ and $C$. Then there exist $C^*\in\mathbb{R}$
such
\begin{equation*}
h(x)=C^*e^{-ax^{2}}.
\end{equation*}
\end{lemma}

The reader can see the reference \cite{Si} for the proof.

\bigskip

Now we are in a position to state and prove the $q$-analogue of the Hardy's
theorem

\begin{theorem}
\label{t2} Suppose $f\in\mathcal{L}_{q,1,\nu}$ satisfying the following
estimates
\begin{equation}  \label{e3}
|f(x)|\leq C e^{-\frac{1}{2} x^2},\quad\forall x\in\mathbb{R}_q^+,
\end{equation}

\begin{equation*}
|\mathcal{F}_{q,\nu} f(x)|\leq C e^{-\frac{1}{2} x^2},\quad\forall x\in%
\mathbb{R},
\end{equation*}
where $C$ is a positive constant. Then there exist $A\in\mathbb{R}$ such
that
\begin{equation*}
f(z)=Ac_{q,\nu}\mathcal{F}_{q,\nu}\left(e^{-\frac{1}{2} x^2}\right)(z),
\quad\forall z\in\mathbb{C}.
\end{equation*}
\end{theorem}

\begin{proof}
We claim that $\mathcal{F}_{q,\nu} f$ is an analytic function and there
exist $C^{\prime}>0$ such that

\begin{equation*}
|\mathcal{F}_{q,\nu }f(z)|\leq C^{\prime }e^{\frac{1}{2}|z|^{2}},\quad
\forall z\in \mathbb{C}.
\end{equation*}
We have
\begin{equation*}
|\mathcal{F}_{q,\nu }f(z)|\leq c_{q,\nu }\int_{0}^{\infty }|f(x)||j_{\nu
}(zx,q^{2})|x^{2\nu +1}d_{q}x.
\end{equation*}
From the Lemma \ref{l1}, if $\left\vert z\right\vert >1$ then there exist $%
\sigma _{1}>0$ such that
\begin{equation*}
x^{2\nu +1}|j_{\nu }(zx,q^{2})|=\frac{1}{\left\vert z\right\vert ^{2\nu +1}}%
(\left\vert z\right\vert x)^{2\nu +1}|j_{\nu }(zx,q^{2})|<\frac{\sigma _{1}}{%
1+\left\vert z\right\vert ^{2}x^{2}}e^{x|z|},\quad \forall x\in \mathbb{R}%
_{q}^{+}.
\end{equation*}
Then we obtain
\begin{equation*}
|\mathcal{F}_{q,\nu }f(z)|\leq C\sigma _{1}c_{q,\nu }\left[ \int_{0}^{\infty
}\frac{e^{-\frac{1}{2}(x-|z|)^{2}}}{1+\left\vert z\right\vert ^{2}x^{2}}%
d_{q}x\right] e^{\frac{1}{2}|z|^{2}}<C\sigma _{1}c_{q,\nu }\left[
\int_{0}^{\infty }\frac{1}{1+x^{2}}d_{q}x\right] e^{\frac{1}{2}|z|^{2}}.
\end{equation*}
Now, if $\left\vert z\right\vert \leq 1$ then there exist $\sigma _{2}>0$
such that%
\begin{equation*}
x^{2\nu +1}|j_{\nu }(zx,q^{2})|\leq \sigma _{2}e^{x},\quad \forall x\in
\mathbb{R}_{q}^{+}.
\end{equation*}
Therefore
\begin{equation*}
|\mathcal{F}_{q,\nu }f(z)|\leq C\sigma _{2}c_{q,\nu }\left[ \int_{0}^{\infty
}e^{-\frac{1}{2}x^{2}+x}d_{q}x\right] \leq C\sigma _{2}c_{q,\nu }\left[
\int_{0}^{\infty }e^{-\frac{1}{2}x^{2}+x}d_{q}x\right] e^{\frac{1}{2}%
|z|^{2}},
\end{equation*}
which leads to the estimate (\ref{e3}). Using Lemma \ref{l2}, we obtain
\begin{equation*}
\mathcal{F}_{q,\nu }f(z)=\text{const}.e^{-\frac{1}{2}z^{2}},\quad \forall
z\in \mathbb{C},
\end{equation*}
and by Theorem \ref{t1}, we conclude that
\begin{equation*}
f(z)=\text{const}.\mathcal{F}_{q,\nu }\left( e^{-\frac{1}{2}t^{2}}\right)
(z),\quad \forall z\in \mathbb{C}.
\end{equation*}
\end{proof}

\begin{corollary}
\label{c1} Suppose $f\in\mathcal{L}_{q,1,\nu}$ satisfying the following
estimates
\begin{equation*}
|f(x)|\leq C e^{-p x^2},\quad\forall x\in\mathbb{R}_q^+,
\end{equation*}

\begin{equation*}
|\mathcal{F}_{q,\nu} f(x)|\leq C e^{-\sigma x^2},\quad\forall x\in\mathbb{R},
\end{equation*}
where $C,p,\sigma$ are a positive constant and $p\sigma=\frac{1}{4}$. We
suppose that there exist $a\in\mathbb{R}_q^+$ such that $a^2p=\frac{1}{2}$.
Then there exist $A\in\mathbb{R}$ such that
\begin{equation*}
f(z)=Ac_{q,\nu}\mathcal{F}_{q,\nu}\left(e^{-\sigma
t^2}\right)(z),\quad\forall z\in\mathbb{C}.
\end{equation*}
\end{corollary}

\begin{proof}
Let $a\in\mathbb{R}_q^+$, and put
\begin{equation*}
f_a(x)=f(ax),
\end{equation*}
then
\begin{equation*}
\mathcal{F}_{q,\nu} f_a(x)=\frac{1}{a^{2\nu+2}}\mathcal{F}_{q,\nu} f(x/a).
\end{equation*}
In the end, applying Theorem \ref{t2} to the function $f_a$.
\end{proof}

\begin{corollary}
Suppose $f\in\mathcal{L}_{q,1,\nu}$ satisfying the following estimates
\begin{equation*}
|f(x)|\leq C e^{-p x^2},\quad\forall x\in\mathbb{R}_q^+,
\end{equation*}

\begin{equation}  \label{e8}
|\mathcal{F}_{q,\nu} f(x)|\leq C e^{-\sigma x^2},\quad\forall x\in\mathbb{R},
\end{equation}
where $C,p,\sigma$ are a positive constant and $p\sigma>\frac{1}{4}$. We
suppose that there exist $a\in\mathbb{R}_q^+$ such that $a^2p=\frac{1}{2}$.
Then $f\equiv 0$.
\end{corollary}

\begin{proof}
In fact there exists $\sigma^{\prime}<\sigma$ such that $p\sigma^{\prime}=
\frac{1}{4}$. Then the function $f$ satisfying the estimates of Corollary %
\ref{c1}, if we replacing $\sigma$ by $\sigma^{\prime}$. Which implies
\begin{equation*}
\mathcal{F}_{q,\nu} f(x)=\text{const}.e^{-\sigma^{\prime}x^2},\quad\forall
x\in\mathbb{R}.
\end{equation*}
On the other hand, $f$ satisfying the estimates (\ref{e8}), then
\begin{equation*}
\left|\text{const}.e^{-\sigma^{\prime}x^2}\right|\leq C e^{-\sigma
x^2},\quad\forall x\in\mathbb{R}.
\end{equation*}
This implies $\mathcal{F}_{q,\nu} f\equiv 0$, and by Theorem \ref{t1} we
conclude that $f\equiv 0$.
\end{proof}

\section{The $q$-Fourier-Neumann Expansions}

The little $q$-Jacobi polynomials are defined for $\nu ,\beta >-1$ by \cite%
{KoS}
\begin{equation*}
p_{n}(x;q^{\nu },q^{\beta };q)={_{2}\phi _{1}}\left( \,\left.
\begin{matrix}
q^{n+\nu +\beta +1},q^{-n} \\
q^{\nu +1}%
\end{matrix}
\right\vert \,q;qx\,\right) .
\end{equation*}
We define the functions
\begin{equation*}
P_{\nu ,n}(x;q^{2})=\sigma _{q,\nu }(n)q^{-n(\nu +1)}\frac{(q^{2+2n},q^{2\nu
+2};q^{2})_{\infty }}{(q^{2+2n+2\nu },q^{2};q^{2})_{\infty }}
p_{n}(x^{2};q^{2\nu },1;q^{2})
\end{equation*}
and
\begin{equation*}
\mathcal{J}_{\nu ,n}(x;q^{2})=\sigma _{q,\nu }(n)\frac{J_{\nu
+2n+1}(q^{n}x;q^{2})}{x^{\nu +1}},
\end{equation*}
where
\begin{equation*}
\sigma _{q,\nu }(n)=\sqrt{\frac{1-q^{2\nu +4n+2}}{1-q}}.
\end{equation*}
Consider $\mathcal{L}_{q,2}^{\nu }$ as an Hilbert space with the inner
product
\begin{equation*}
\langle f|g\rangle =\int_{0}^{1}f(x)g(x)x^{2\nu+1}d_{q}x.
\end{equation*}
The $q$-Paley-Wiener space is defined by
\begin{equation*}
PW_{q}^{\nu}=\left\{ f\in \mathcal{L}_{q,2,\nu }:f(x)=c_{q,\nu
}\int_{0}^{1}u(t)j_{\nu }(xt,q^{2})t^{2\nu +1}d_{q}t,\quad u\in \mathcal{L}%
_{q,2}^{\nu }\right\} .
\end{equation*}

\begin{proposition}
$PW_{q}^{\nu}$ is a closed subspace of $\mathcal{L}_{q,2,\nu }$ and with the
inner product given in (\ref{e4}) is an Hilbert space.
\end{proposition}

\begin{proof}
In fact, given $f\in \mathcal{L}_{q,2,\nu }$ and let $\{f_{n}\}_{n\in
\mathbb{N}}$ be a sequence of element of $PW_{q}^{\nu}$ which converge to $f$
in $L^{2}$-norm. For $n\in \mathbb{N}$, there exist $u_{n}\in \mathcal{L}
_{q,2}^{\nu }$ such that
\begin{equation*}
f_{n}(x)=c_{q,\nu }\int_{0}^{1}u_{n}(t)j_{\nu }(xt,q^{2})t^{2\nu +1}d_{q}t.
\end{equation*}
Moreover
\begin{equation*}
\lim_{n\rightarrow \infty }\Vert f_{n}-f\Vert _{q,2,\nu }=0.
\end{equation*}
This give
\begin{equation*}
\lim_{n\rightarrow \infty }\Vert \mathcal{F}_{q,\nu }f_{n}-\mathcal{F}%
_{q,\nu }f\Vert _{q,2,\nu }=0,
\end{equation*}
and then
\begin{equation*}
\lim_{n\rightarrow \infty }\left[ \int_{0}^{1}|\mathcal{F}_{q,\nu }f_{n}(x)-%
\mathcal{F}_{q,\nu }f(x)|^{2}x^{2\nu +1}d_{q}x+\int_{1}^{\infty }|\mathcal{F}%
_{q,\nu }f(x)|^{2}x^{2\nu +1}d_{q}x\right] =0,
\end{equation*}
which implies
\begin{equation*}
\int_{1}^{\infty }|\mathcal{F}_{q,\nu }f(x)|^{2}x^{2\nu
+1}d_{q}x=0\Rightarrow \mathcal{F}_{q,\nu }f(x)=0,\quad \forall x\in \mathbb{R}_q^+\cap ]1,+\infty[.
\end{equation*}
Then $f\in PW_{q}^{\nu}$.
\end{proof}

\begin{proposition}
We have
\begin{equation*}
\mathcal{F}_{q,\nu }(\mathcal{J}_{\nu ,n})(x)=P_{\nu ,n}(x;q^{2})\chi
_{\lbrack 0,1]}(x),\quad \forall x\in \mathbb{R}_q^+.
\end{equation*}
As a consequence
\begin{equation*}
\int_{0}^{1}P_{\nu ,n}(x;q^{2})P_{\nu ,m}(x;q^{2})x^{2\nu
+1}\,d_{q}x=\delta_{n,m}.
\end{equation*}
\end{proposition}

\begin{proof}
The following proof is identical to the proof of Lemma 1 in \cite{A}. Using
an identity established in \cite{K,K1}
\begin{multline}  \label{e5}
\int_{0}^{\infty }t^{-\lambda }J_{\mu}(q^{m}t;q^{2})J_{\theta }(q^{n}t;q^{2})d_{q}t \\
=(1-q)q^{n(\lambda -1)+(m-n)\mu }\frac{(q^{1+\lambda +\theta -\mu },q^{2\mu
+2};q^{2})_{\infty }}{(q^{1-\lambda +\theta +\mu },q^{2};q^{2})_{\infty }} \\
\times {_{2}\phi _{1}}\left( \,\left.
\begin{matrix}
q^{1-\lambda +\mu +\theta },q^{1-\lambda +\mu -\theta } \\
q^{2\mu +2}%
\end{matrix}%
\right\vert \,q^{2};q^{2m-2n+1+\lambda +\theta -\mu }\,\right) ,
\end{multline}
where $n,m\in \mathbb{Z}$ and $\theta ,\mu ,\lambda \in \mathbb{C}$ such
that ${\rm Re}(1-\lambda +\theta +\mu) >0$, $\theta, \mu$ are not equal to a negative integer and
$$
(\lambda+\theta+1-\mu)/2,\quad m-n+(\lambda+\theta+1-\mu)/2
$$
are not a non-positive integer \cite{K1}.\bigskip

To evaluate $\mathcal{F}_{q,\nu }(\mathcal{J}_{\nu ,n})(x)$ when $%
x=q^{m}\leq 1$, we take in (\ref{e5})
\begin{equation*}
q^{m}=x,\mu =\nu ,\theta =\nu +2n+1,\lambda =0
\end{equation*}%
then
\begin{eqnarray*}
\mathcal{F}_{q,\nu }(\mathcal{J}_{\nu ,n})(x) &=&\sigma _{q,\nu }(n)\frac{%
x^{-\nu }}{1-q}\int_{0}^{\infty }J_{\nu }(xt;q^{2})J_{\nu
+2n+1}(q^{n}t;q^{2})\,d_{q}t \\
&=&\sigma _{q,\nu }(n)q^{-n(\nu +1)}\frac{(q^{2+2n},q^{2\nu
+2};q^{2})_{\infty }}{(q^{2+2n+2\nu },q^{2};q^{2})_{\infty }}{\,{}_{2}\phi
_{1}}\left( \,\left.
\begin{matrix}
q^{2+2\nu +2n},q^{-2n} \\
q^{2\nu +2}%
\end{matrix}%
\right\vert \,q^{2};q^{2}x^{2}\,\right) \\
&=&P_{\nu ,n}(x;q^{2}).
\end{eqnarray*}%
To evaluate $\mathcal{F}_{q,\nu }(\mathcal{J}_{\nu ,m})(x)$ when $x=q^{n}>1$%
, we consider in (\ref{e5})
\begin{equation*}
q^{n}=x,\mu =\nu +2m+1,\theta =\nu ,\lambda =0
\end{equation*}%
In this way, $1+\lambda +\theta -\mu =-2m$. This gives for $m\in \mathbb{N}$
a factor
\begin{equation*}
(q^{-2m};q^{2})_{\infty }=0
\end{equation*}%
on the numerator and then
\begin{equation*}
\mathcal{F}_{q,\nu }(\mathcal{J}_{\nu ,m})(x)=0,\quad x>1
\end{equation*}%
By setting $\lambda =1$, $\theta =\nu +2n+1$, and $\mu =\nu +2m+1$ in~, it
is clear that, for $n,m=0,1,2,\dots $,
\begin{equation*}
\int_{0}^{\infty }J_{\nu +2n+1}(q^{n}x;q^{2})J_{\nu +2m+1}(q^{m}x;q^{2})\,%
\frac{d_{q}x}{x}=\frac{1}{\sigma _{q,\nu }(n)^{2}}\delta _{n,m}
\end{equation*}%
and then
\begin{equation*}
\int_{0}^{\infty }\mathcal{J}_{\nu ,n}(x;q^{2})\mathcal{J}_{\nu
,m}(x;q^{2})\,x^{2\nu +1}d_{q}x=\delta _{n,m}.
\end{equation*}%
Now we use the arguments of $q$-Bessel Fourier analysis provided in this paper
to show that
\begin{equation}
\left\langle P_{\nu ,n}\chi _{\lbrack 0,1]},P_{\nu ,m}\chi _{\lbrack
0,1]}\right\rangle =\left\langle \mathcal{F}_{q,\nu }(\mathcal{J}_{v,n}),%
\mathcal{F}_{q,\nu }(\mathcal{J}_{v,m})\right\rangle =\left\langle \mathcal{J%
}_{\nu ,n},\mathcal{J}_{\nu ,m}\right\rangle =\delta _{n,m}.  \label{e6}
\end{equation}
Another proof of the orthogonality of the little $q$-Jacobi polynomials can
be found in \cite{KoS}
\end{proof}

\begin{proposition}
The systems
\begin{equation*}
\{\mathcal{J}_{\nu,n}\}_{n=0}^\infty,\quad \{P_{\nu,n}\}_{n=0}^\infty
\end{equation*}
forme two orthonormals basis respectively of the Hilbert spaces $PW^v_{q}$
and $\mathcal{L}_{q,2}^\nu$.
\end{proposition}

\begin{proof}
From (\ref{e6}) we derive the orthonormality. To prove that the system $\{%
\mathcal{J}_{\nu ,n}\}_{n=0}^{\infty }$ is complet in $PW_{q}^{\nu}$, given
a function $f\in PW_{q}^{\nu}$ such that
\begin{equation*}
\left\langle f,\mathcal{J}_{\nu ,n}\right\rangle =0,\quad \forall n\in
\mathbb{N}.
\end{equation*}%
Then
\begin{equation*}
\left\langle \mathcal{F}_{q,\nu }(f),\mathcal{F}_{q,\nu }(\mathcal{J}_{\nu
,n}\right)\rangle=0,\quad \forall n\in \mathbb{N},
\end{equation*}%
which implies
\begin{equation*}
\left\langle \mathcal{F}_{q,\nu }(f),P_{\nu ,n}\chi _{\lbrack
0,1]}\right\rangle =\left\langle \mathcal{F}_{q,\nu }(f)\chi _{\lbrack
0,1]},P_{\nu ,n}\right\rangle =\left\langle \mathcal{F}_{q,\nu }(f),P_{\nu
,n}\right\rangle =0,\quad \forall n\in \mathbb{N}.
\end{equation*}%
From the definition of the polynomial $P_{\nu ,n}$ we conclude that
\begin{equation*}
\left\langle \mathcal{F}_{q,\nu }(f),t^{2n}\right\rangle =0,\quad \forall
n\in \mathbb{N}.
\end{equation*}%
Then
\begin{equation*}
c_{q,\nu }\sum_{n=0}^{\infty }(-1)^{n}\frac{q^{n(n+1)}}{%
(q^{2},q^{2})_{n}(q^{2\nu +2},q^{2})_{n}}\left\langle \mathcal{F}_{q,\nu
}(f),t^{2n}\right\rangle x^{2n}=0,\quad \forall x\in \mathbb{R}_{q}^{+},
\end{equation*}%
which can be written as
\begin{equation*}
\mathcal{F}_{q,\nu }^{2}(f)(x)=0,\quad \forall x\in \mathbb{R}_{q}^{+}.
\end{equation*}%
By the inversion formula (\ref{e1}) we conclude that $f=0$. From (\ref{e6})
we derive the orthonormality. To prove that the system $\{P_{\nu
,n}\}_{n=0}^{\infty }$ is complet in $\mathcal{L}_{q,2}^{\nu }$, given a
function $f\in \mathcal{L}_{q,2}^{\nu }$ such that
\begin{equation*}
\langle f|P_{\nu ,n}\rangle =0,\quad \forall n\in \mathbb{N}
\end{equation*}%
Then
\begin{equation*}
\langle f|t^{2n}\rangle =0,\quad \forall n\in \mathbb{N}.
\end{equation*}
Which leads to the result.
\end{proof}

\begin{proposition}
Let $\lambda \in \mathbb{R}_{q}^{+}$ then
\begin{equation*}
c_{q,\nu }j_{\nu }(\lambda x;q^{2})=\sum_{n=0}^{\infty }\mathcal{J}_{n,\nu
}(\lambda ;q^{2})P_{n,\nu }(x),\quad \forall x\in \lbrack 0,1]\cap \mathbb{R}%
_{q}^{+}.
\end{equation*}%
As a consequence we have
\begin{equation*}
\sum_{n=0}^{\infty }\left[ P_{n,\nu }(x;q^{2})\right] ^{2}=\frac{x^{-2(\nu
+1)}}{1-q},\quad \forall x\in \lbrack 0,1]\cap \mathbb{R}_{q}^{+}
\end{equation*}%
and for all $\lambda \in \mathbb{R}_{q}^{+}$
\begin{eqnarray*}
&&\sum_{n=0}^{\infty }\left[ \mathcal{J}_{n,\nu }(\lambda ;q^{2})\right]
^{2}=-\frac{q^{\nu }}{2(1-q)\lambda ^{1+2\nu }} \\
&&\times \left[ \frac{\lambda }{q}J_{\nu +1}(\lambda ;q^{2})J_{\nu }^{\prime
}(\lambda /q;q^{2})-J_{\nu +1}(\lambda ;q^{2})J_{\nu }(\lambda
/q;q^{2})-J_{\nu +1}^{\prime }(\lambda ;q^{2})J_{\nu }(\lambda /q;q^{2})%
\right].
\end{eqnarray*}
\end{proposition}

\begin{proof}
Let $\lambda \in \mathbb{R}$ and consider the function
\begin{equation*}
\psi _{\lambda }:[0,1]\cap \mathbb{R}_{q}^{+}\rightarrow \mathbb{R},\quad
x\mapsto c_{q,\nu }j_{\nu }(\lambda x;q^{2}).
\end{equation*}
Then $\psi _{\lambda }\in \mathcal{L}_{q,2}^{\nu }$ and we can write
\begin{equation}  \label{e7}
\psi _{\lambda }(x)=\sum_{n=0}^{\infty }\left\langle \psi _{\lambda
}|P_{n,\nu }\right\rangle P_{n,\nu }(x),\quad \forall x\in \lbrack 0,1]\cap
\mathbb{R}_{q}^{+}.
\end{equation}
Note that
\begin{equation*}
\left\langle \psi _{\lambda }|P_{n,\nu }\right\rangle =\left\langle \psi
_{\lambda },P_{n,\nu }\chi _{\lbrack 0,1]}\right\rangle =\left\langle \psi
_{\lambda },\mathcal{F}_{q,\nu }(\mathcal{J}_{n,\nu })\right\rangle =%
\mathcal{F}_{q,\nu }^{2}(\mathcal{J}_{n,\nu })(\lambda )=\mathcal{J}_{n,\nu
}(\lambda ;q^{2}).
\end{equation*}%
Then we deduce the result. Using the Parseval's theorem and (\ref{e7}) we
obtain
\begin{equation*}
\sum_{n=0}^{\infty }\left[ P_{n,\nu }(x;q^{2})\right] ^{2}=\Vert \psi
_{x}\Vert _{q,2,\nu }^{2}=\frac{x^{-2(\nu +1)}}{1-q}.
\end{equation*}%
The second identity is deduced also from the Parseval's theorem
\begin{equation*}
\sum_{n=0}^{\infty }\left[ \mathcal{J}_{n,\nu }(\lambda ;q^{2})\right]
^{2}=N_{q,\nu ,2}^{2}(\psi _{\lambda }),
\end{equation*}%
and the following relation proved in \cite{Ko}
\begin{eqnarray*}
&&\int_{0}^{1}\left[ J_{\nu }(aqt;q^{2})\right] ^{2}td_{q}t=-\frac{%
(1-q)q^{\nu -1}}{2a} \\
&&\times \Big[ aJ_{\nu +1}(aq;q^{2})J_{\nu }^{\prime }(a;q^{2})-J_{\nu
+1}(aq;q^{2})J_{\nu }(a;q^{2})-J_{\nu +1}^{\prime }(aq;q^{2})J_{\nu
}(a;q^{2})\Big] .
\end{eqnarray*}
\end{proof}

\end{document}